# METHOD OF MONOTONE STRUCTURAL EVOLUTION FOR CONTROL AND STATE CONSTRAINED OPTIMAL CONTROL PROBLEMS

**Maciej Szymkat, Adam Korytowski**

AGH University of Science and Technology
30-059 Kraków, Poland, e-mail: msz@ia.agh.edu.pl

## Abstract

A method of optimal control computation is proposed for problems with control and state constraints. It uses a sequence of control structure adjustments in the form of generations and reductions of nodes and arcs, which do not change the current control but redefine the decision space. Several examples are given.

## 1 Introduction

Numerical methods of optimal control are divided into *direct* and *indirect* [1]. In the direct approach an approximating, finite dimensional optimization problem is constructed and solved by nonlinear programming algorithms. Collocation and SQP methods are often used, usually leading to large-scale computations. The well established direct methods feature large areas of convergence but they are rather slow, especially in the final stage. This approach is not particularly demanding upon the user and is considered fairly universal. It has many powerful implementations like SOCS (Betts [1]), DIRCOL (von Stryk [10]) or DIRMUS (Hinsberger).

In the indirect approach the optimal solution is computed by solving the boundary value problem obtained from the maximum principle. Multiple shooting is frequently used, with such implementations as BNDSCO (Oberle and Grimm), MUMUS (Hiltmann *et al.*) and MUSCOD-II (Diehl). The collocation methods for indirect computations (Kierzenka and Shampine [4]) involve large systems of algebraic equations requiring specialized algorithms. The rate of convergence of indirect methods is usually very high, but their area of convergence is small and so they require good initial guesses for the adjoint vector. Practically, the optimal control structure has to be known beforehand. This can be achieved by a direct algorithm (Shen and Tsiotras [8]) or using *homotopy* methods where a sequence of appropriately constructed auxiliary problems is solved by multiple shooting (Bulirsch *et al.* [3]).

This paper describes a new direct approach to numerical dynamic optimization, called *monotone structural evolution* (MSE). It is effective for a large class of nonlinear problems with control and state constraints, including singular cases. MSE originates from the *variable parameterization* method, developed for the nonsingular affine case with simple control bounds in [5, 6, 11 – 15, 17]. A distinctive feature of MSE is that the decision space is systematically reconstructed in the course of optimization, with changing the control structure, parameterization and, typically, the number of decision variables. The search for structural changes which lead to rapid improvement of the performance index is based on analysis of the discrepancy between the current approximation of solution and the maximum principle optimality conditions, and continues until these conditions are satisfied with sufficient accuracy. The proper choice of the sequence of decision spaces, utilizing information taken from the adjoint solution, allows the number of decision variables to be kept comparatively small, at least in early stages of optimization. In consequence, quasi-Newton or Newton optimization with analytical gradients can be used, which is vital for fast convergence. The dimension of the decision space grows only when this is necessary for improving accuracy of optimal control approximation. An important property of MSE is that the performance index decreases monotonously during optimization, due to control preservation by the structure changes.

## 2 Optimal control problem

The control system is described by a state equation

$$\dot{x} = f(x,u), \quad t \in [0,T], \quad x(0) = x^0 \qquad (2.1)$$

where the state $x(t) \in \mathbf{R}^n$ and control $u(t) \in \mathbf{R}^m$. The admissible controls are piecewise continuous functions satisfying a vector inequality

$$g(x,u) \le 0. \qquad (2.2)$$

It is assumed that for each component of $g$, the equation $g_i(x,u) = 0$ always has a unique solution w.r.t. some component of $u$. The optimal control problem is to find an admissible control minimizing a performance functional

$$S(u) = \varphi(x(T)). \qquad (2.3)$$

The initial state $x^0$ and the horizon $T$ are fixed. Vector state constraints of the form $h(x) \le 0$ and terminal conditions on state can be also introduced and treated by means of penalty techniques. The functions $f$, $g$, $h$ and $\varphi$ are sufficiently

smooth. The adjoint trajectory $\psi$ is defined as a solution of the boundary value problem

$$\dot{\psi} = -\nabla_x H(\psi, x, u), \quad \psi(T) = -\nabla \varphi(x(T)) \qquad (2.4)$$

where $H(\psi, x, u) = \psi^\top f(x, u)$ is the hamiltonian.

## 3 Monotone Structural Evolution

### 3.1 Control structure

The structure of optimal control is usually understood as the sequence of sets of active constraints. The elements in the sequence are arranged as they follow in time. The corresponding segments of control and state trajectory are called arcs (boundary, constrained, singular, etc.). We introduce a more general concept of structure for arbitrary admissible controls, basing on the observation that active constraints determine the shape of optimal control and the way it is computed. The *control structure* will thus mean the sequence of procedures used to determine control values in successive time intervals. The procedures for control computation are built with the condition of hamiltonian maximization taken into account. Typically, the control structure varies during optimization.

For every control with a given structure a sequence of *division nodes* $\tau_1, \ldots, \tau_N$ is defined,

$$0 = \tau_0 \leq \tau_1 \leq \tau_2 \leq \ldots \leq \tau_N = T. \qquad (3.1)$$

In each interval $[\tau_{i-1}, \tau_i[$, a fixed procedure $P_i$ is used to calculate control. The restrictions of control and state trajectory to intervals $[\tau_{i-1}, \tau_i[$ are called *arcs*. The *parameters of a control structure* include the division nodes and, possibly, other parameters used by the procedures $P_i$. The structure parameters are decision variables in the optimization algorithms of MSE. As quasi-Newton or Newton gradient optimization algorithms are used in the sequel, we need analytical formulas for derivatives of the performance functional (2.3) with respect to decision variables. It is convenient to express them in terms of the adjoint trajectory.

### 3.2 Generations and reductions

An optimal control approximation in the direct approach is a value of a mapping $P: D \rightarrow U$ from a finite-dimensional space of decision variables into the functional control space. Every structural change in MSE, that is, both the *generation* and *reduction* basically consists in choosing a new decision space $D'$ and a new function $P': D' \rightarrow U$. It is required that the control is not immediately affected, which means that $P'(d') = P(d)$ for $d \in D$ and $d' \in D'$ being the decision vectors just before and after the structural change. Thanks to this *condition of control preservation*, the performance index monotonously decreases during the overall optimization procedure. The dimension of the decision space increases in a generation, and is diminished in a reduction. Typically, only few selected elements defining the structure can be affected by a structural change. For example, one or two new procedures $P'_i$ are introduced in a generation with inserting the corresponding new nodes, or one of the procedures $P_i$ is modified.

A typical reduction consists of eliminating an arc of zero length. More precisely, every arc of zero length is subject to reduction if the directional derivative of performance index w.r.t. its boundaries is nonnegative for all admissible directions. At the same time the decision variables that describe this control arc are eliminated, including at least one of the respective division nodes. Such a reduction occurs each time when one of the constraints (3.1) becomes active after the linesearch of the gradient optimization algorithm. Another typical reduction occurs when two adjacent arcs described by identical procedures are unified.

### 3.3 General algorithm

The basic algorithm of MSE consists of the following steps.

$1^0$ Selection of starting point.

$2^0$ Termination, if optimality conditions are satisfied.

$3^0$ Generation, if it is sufficiently promising or needed.

$4^0$ Iteration of gradient optimization.

$5^0$ Reduction, if necessary.

$6^0$ Return to $2^0$.

Step $3^0$ is distinctive for MSE algorithms and crucial for their convergence. The changes of structure are mainly performed to speed up the optimization when a stationary point in the current decision space is being approached. This algorithm should be equipped with special procedures for gradient computation and evaluation of efficiency of generations. To treat state constraints, an outer loop of penalty modification has to be added. In the gradient optimization of step $4^0$ the bounds on division nodes and control constraints are respected due to an appropriate organization of linesearch. Numerical solutions of differential equations are obtained by the RK4 method with mesh adjusted so as to include all discontinuity points.

### 3.4 Rules for generations

The *efficiency* of a generation is defined as the difference of squared Euclidean norms of antigradients of the performance index w.r.t. the decision vector immediately after and before the generation. This definition applies only to the case of admissible antigradients. In the general case the antigradients are replaced by their orthogonal projections onto the local conical approximation of the admissible set. Such a definition is justified in two ways. First, the square of the norm of gradient, multiplied by $-1$, is equal to the derivative of the performance index w.r.t. the search line parameter in the steepest descent direction. The efficiency thus determines the increase of steepness of the performance index in this direction. Secondly, the efficiency so defined does not depend on those components of the gradient of performance index that are not affected by the generation, which simplifies computations.

There are two kinds of generations: intentional, aimed at speeding up the optimization, and those enforced by the requirement that at the moment of gradient computation each control arc has to be either purely boundary or purely interior. The latter, called *saturation generations*, are performed when the optimization transforms an interior division interval into one that contains a segment with active control constraint (2.2). Such an interval is divided by introducing new division nodes.

The generations of the first kind (intentional) are primarily based on relative efficiency, that is, the ratio of efficiency and the squared gradient norm. The generation takes place if the relative efficiency exceeds a given threshold. By choosing the threshold one can control the trade-off between the dimension of decision space and gradient magnitude. The number of simultaneously generated nodes is limited by additional rules (e.g, one or two per arc, solely at local maximizers of relative efficiency), to avoid an undesirable increase of the number of decision variables. Additional requirements can be imposed on generations to obtain controls with pre-selected regularity properties, like continuity or smoothness. The choice of particular generations is also subject to the condition that optimization should converge to the optimal control in the strong sense.

## 4 Examples

### 4.1 LQ problem with control constraints

Consider a simple control-constrained LQ problem

$\dot{x}_1 = x_2$, $\dot{x}_2 = -x_1 + u$, $t \in [0,T]$

$x(0) = \text{col}(4, -4)$, $|u(t)| \leq 1$

$S(u) = \frac{1}{2} \int_0^T (x_1^2 + x_2^2 + u^2) dt$, $T = 15$.

The hamiltonian $H = \psi_1 x_2 + \psi_2(-x_1 + u) - \frac{1}{2}(x_1^2 + x_2^2 + u^2)$ attains maximum at

$u = \begin{cases} \psi_2, & |\psi_2| \leq 1 \\ \text{sgn}\,\psi_2, & |\psi_2| > 1 \end{cases}$

where $\psi$ satisfies $\dot{\psi}_1 = x_1 + \psi_2$, $\dot{\psi}_2 = x_2 - \psi_1$, $\psi(T) = 0$. The optimal control is continuous with two types of arcs possible: boundary and interior. Its first derivative may be discontinuous only at the ends of boundary arcs. It is assumed that admissible approximations of optimal control also have these properties. In every interior arc

$u(t) = p_i^T w(t, \tau_{i-1}, \tau_i)$ (4.1)

where $p_i$ is a vector of parameters and $w$ is a given vector function. First, we take the Hermite cubic polynomials

$w_1(t, \tau_{i-1}, \tau_i) = w_3(t, \tau_i, \tau_{i-1}) = (t-\tau_i)^2(2t+\tau_i - 3\tau_{i-1})/(\tau_i - \tau_{i-1})^3$

$w_2(t, \tau_{i-1}, \tau_i) = w_4(t, \tau_i, \tau_{i-1}) = (t-\tau_i)^2(t-\tau_{i-1})/(\tau_i - \tau_{i-1})^2$,

thus $u(\tau_{i-1}) = p_{i1}$, $\dot{u}(\tau_{i-1}+) = p_{i2}$, $u(\tau_i) = p_{i3}$, $\dot{u}(\tau_i-) = p_{i4}$.
To ensure continuity at division points and smoothness between neighboring interior arcs, some parameters $p_{ik}$ are fixed or identified.

Let $\Sigma$ denote the performance index as a function of the parameters and division points. Its derivative w.r.t. $p_{ik}$ reads

$$\nabla_{p_{ik}} \Sigma = -\int_{\Pi_{ik}} \nabla_u H(\psi, x, u) \nabla_{p_{ik}} u \, dt \qquad (4.2)$$

where the derivatives of $u$ are determined by (4.1), and $\Pi_{ik}$ is the union of $[\tau_{i-1}, \tau_i]$ and, possibly, one of its neighboring interior intervals. The derivative w.r.t. $\tau_i$, $i = 1,...,N-1$ being the right-hand end of an interior interval is given by

$\nabla_{\tau_i} \Sigma = -\dot{u}(\tau_i-)\nabla_{p_{i3}}\Sigma - \ddot{u}(\tau_i-)\nabla_{p_{i4}}^-\Sigma - \ddot{u}(\tau_i+)\nabla_{p_{i4}}^+\Sigma$

where $\nabla_{p_{i4}}^-\Sigma$ and $\nabla_{p_{i4}}^+\Sigma$ are computed according to (4.2), but with $\Pi_{ik}$ equal to $[\tau_{i-1}, \tau_i]$ and $[\tau_i, \tau_{i+1}]$, respectively. For the left-hand end points we have

$\nabla_{\tau_i} \Sigma = -\dot{u}(\tau_i+)\nabla_{p_{i+1,1}}\Sigma - \ddot{u}(\tau_i-)\nabla_{p_{i+1,2}}^-\Sigma - \ddot{u}(\tau_i+)\nabla_{p_{i+1,2}}^+\Sigma$.

If $\tau_i$ is an end point of a boundary arc, the terms with vanishing control derivatives are dropped.

The optimization starts from zero control. If we only admit generations of arcs due to control saturation, we obtain the approximation of optimal control shown with solid line in the upper part of Fig. 1. The division points are marked with circles. Since $\nabla_u H = \psi_2 - u$, this control is evidently far from optimal.

In the second stage of optimization, generations of additional division points on the interior arcs are allowed. The generation criteria are based on relative efficiency and the distance of the control zeroing $\nabla_u H$ (i.e., equal to $\psi_2$) from its local cubic approximation. The result, close to optimal, is shown in the lower part of Fig. 1. Observe the significant increase of the dimension of decision space (to 60).

The second stage of optimization has been repeated with generations of another type consisting in an extension of the vector $w$ with four new components

$w_{k+4} = v_k - v_k(\tau_{i-1})w_1 - v_k(\tau_i)w_3$, $k = 1,...,4$

where $v_1 = e^{\alpha s} \sin \beta s$, $v_2 = e^{\alpha s} \cos \beta s$, $v_3 = e^{-\alpha s} \sin \beta s$, $v_4 = e^{-\alpha s} \cos \beta s$, $s = t - \tau_{i-1}$ are basis solutions of the canonical system restricted to an interior arc. Note that $u(\tau_{i-1}) = p_{i1}$, $u(\tau_i) = p_{i3}$. The continuity requirement, still valid, is thus easy to impose. The derivatives of $\Sigma$ are given by (4.2) and

$$\nabla_{\tau_i} \Sigma = -\int_{\tau_{i-1}}^{\tau_{i+1}} \nabla_u H(\psi, x, u) \nabla_{\tau_i} u \, dt, \quad i = 1,...,N-1$$

where the derivative of $u$ is by definition zero on boundary intervals.

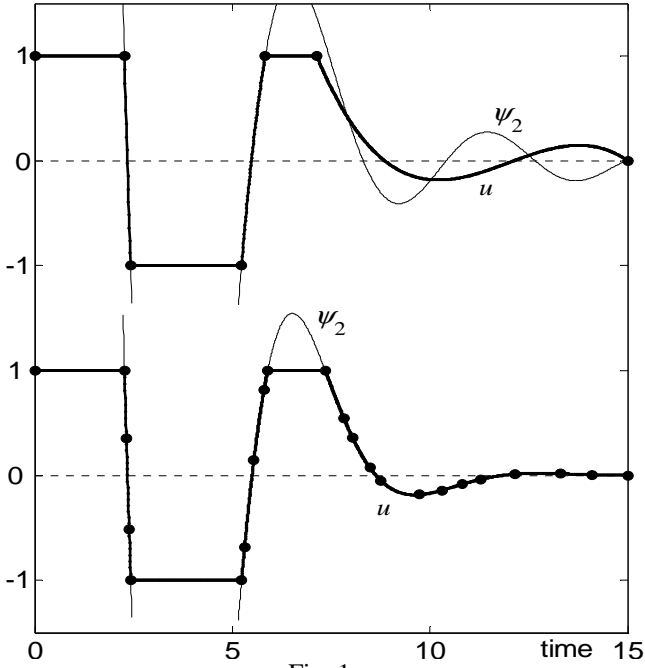

Fig. 1.

The optimization is restarted with such a generation applied only to the last arc. The parameters $p_{ik}$, $k = 1,...,4$ are inherited from the first stage, and $p_{ik} = 0$ for $k = 5,...,8$. The dimension of decision space increases to 15. This stage results in an approximation of the optimal control nearly identical with that shown in the lower part of Fig. 1 (but with only 6 nodes), and the coefficients of the polynomials close to zero for the last arc. Note that this particular choice of approximating functions leads to a solution satisfying the maximum principle conditions without any additional generations during the second stage.

### 4.2 Singular problem

The second example concerns a fermentation model (based on [7]) with a singular arc in the optimal solution. The problem is described as follows

$$\dot{x}_1 = \frac{a_1 x_1 x_2}{1+b_1 x_2 + b_2 x_2^2}\left(\frac{1}{x_3}+c_1+c_2 x_2\right) - \frac{x_1}{x_3} u$$

$$\dot{x}_2 = -\frac{a_2 x_1 x_2}{1+b_3 x_2 + b_4 x_2^2} + \frac{200 - x_2}{x_3} u, \quad \dot{x}_3 = u$$

with $S(u) = (x_2(T) - 50)x_3(T)$, $T = 6$, $x(0) = \text{col}(3, 40, 5)$, $0 \le u(t) \le 1$, $a_1 = 0.2$, $a_2 = 0.5$, $b_1 = 0.2$, $b_2 = 0.001$, $b_3 = 0.1$, $b_4 = 0.0004$, $c_1 = 0.25$, $c_2 = 0.00125$.

As the rhs of the state equations are affine in control, we expect only two kinds of arcs in the optimal solution: boundary and singular. We therefore admit arcs that are either boundary or interior where control is determined from equating the switching function $\phi(x,\psi)$, and its first and second derivatives to zero (singularity of second order). We obtain a formula $u(t) = u_{\text{int}}(x(t))$ expressing this control only by the state variable (which is rather particular). Note that the adjoint equation for the interior arcs reads $\dot{\psi} = -\nabla_x H(\psi, x, u_{\text{int}}(x))$ where $\nabla_x$ denotes the total derivative.

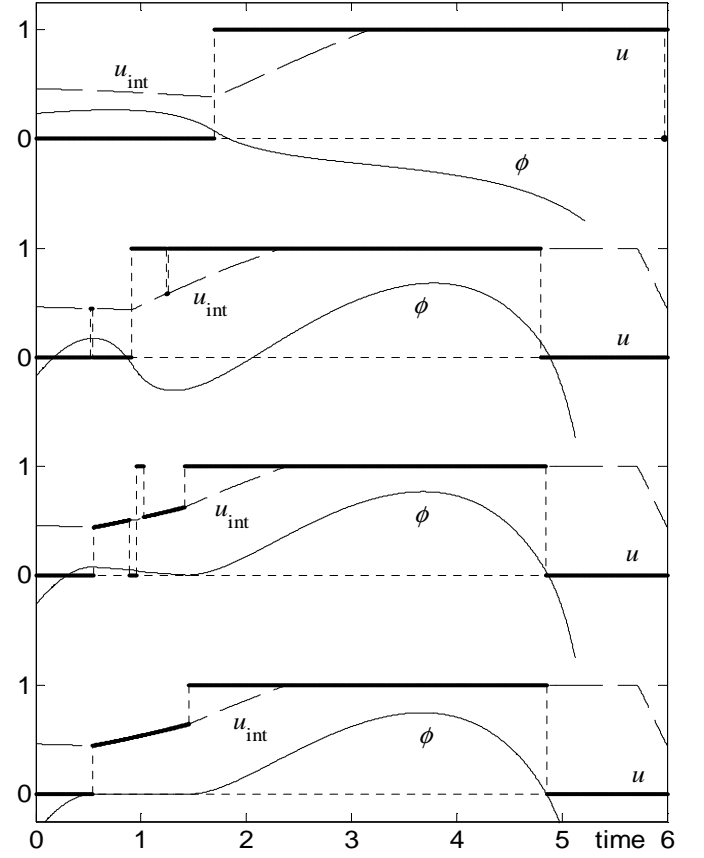

Fig. 2.

In this example, only *spike* and saturation generations are used. A spike generation consists in adding two nodes at the same point $\tau \in \,]0,T[\,$, or one node at $\tau = 0$ or $\tau = T$. A seed of a new arc is thus planted into the structure. No requirement of control continuity on the boundaries of the new arc is imposed. The procedure of control calculation in the new spike is selected as follows. The arc is interior, if

$\phi(x(\tau),\psi(\tau)) > 0$ and $u(\tau) < u_{\text{int}}(x(\tau)) < 1$ or
$\phi(x(\tau),\psi(\tau)) < 0$ and $0 < u_{\text{int}}(x(\tau)) < u(\tau)$.

If neither of these occurs, the upper control bound is taken when $\phi(x(\tau),\psi(\tau)) > 0$ and $u(\tau) < 1$, and the lower bound when $u(\tau) > 0$ and $\phi(x(\tau),\psi(\tau)) < 0$. Besides, the rules of section 3.4 are applied. Effects of weak convergence (chattering control) are eliminated by imposing a lower threshold for the efficiency of generations on interior arcs. The nodes coinciding with discontinuity points are the only decision variables. The derivative $\nabla_{\tau_i}\Sigma$ is equal to the jump of the hamiltonian [9, 16]

$$\Delta H_i = H(\psi(\tau_i), x(\tau_i), u(\tau_i+)) - H(\psi(\tau_i), x(\tau_i), u(\tau_i-)).$$

The optimization is started from a control that switches from 0 to 1 at $t = 3$. During the optimization the control structure evolves in a number of generations and reductions. The first spike generation from the upper to lower bound is shown on the right of the first plot of Fig. 2. After a few iterations in a constant decision space, two other generations to $u_{int}$ occur (second plot). Further optimization leads to an expansion of the new interior arcs (third plot). Finally, after two reductions and combining two interior arcs into one, we obtain the optimal control shown in the last plot of Fig. 2.

### 4.3 State constrained problem

The model below is taken from [2, example 9.3.13]. The task is to swing up an inverted pendulum on a cart from the down stable position, by applying a bounded horizontal force to the cart. The state equations are

$$\dot{x}_1 = x_2, \quad \dot{x}_2 = -(\varepsilon x_2^2 sc + s + cu)/(1 - \varepsilon c^2)$$
$$\dot{x}_3 = x_4, \quad \dot{x}_4 = (\varepsilon s(x_2^2 + c) + u)/(1 - \varepsilon c^2)$$
$$c = \cos x_1, \quad s = \sin x_1, \quad \varepsilon = 0.5.$$

The variable $x_1$ denotes the angle of the pendulum, $x_2$ its angular velocity, $x_3$ the position of the cart and $x_4$, velocity of the cart. The initial state is $x(0) = 0$. The force $u$ and the cart position are bounded, $|u(t)| \leq 1$ and $x_3(t) \leq x_{3\max}$. The performance index has the form

$$S(u) = \tfrac{1}{2}\left((x_1(T) - \pi)^2 + x_2(T)^2 + x_3(T)^2 + x_4(T)^2\right)$$

with $T = 4$.

We use the standard penalty technique by introducing

$$\dot{x}_5 = \begin{cases} 0, & x_3 \leq x_{3\max} \\ \tfrac{1}{2}(x_3 - x_{3\max})^2, & x_3 > x_{3\max} \end{cases}, \quad x_5(0) = 0$$

$$S_\rho(u) = S(u) + \rho x_5(T).$$

As the rhs of the state equations are affine in $u$, only boundary, singular, and state-constrained arcs are possible. By equating the derivatives of $x_3$ and $x_4$ to zero we obtain $u_{con}(x) = -\varepsilon s(x_2^2 + c)$ for candidate state-constrained arcs. The singular arcs are explicitly described in [14]. Since the expression for control depends on the adjoint variables, the approach presented in section 4.2 is not directly applicable. Here we approximate the singular control by explicit functions of time $u_{sng}(t)$ using the cubic Hermite polynomials of section 4.1, but without the continuity or smoothness requirements. The adjoint equation has the form (2.4) with the substitution of $u_{sng}$ and $u_{con}(x)$ for $u$ on respective intervals, and $\nabla_x$ understood as the operator of total differentiation. The derivatives of $\Sigma$ w.r.t. the parameters of polynomial approximation are computed by (4.2), the derivatives w.r.t. the ends of the candidate singular arcs read

$$\nabla_{\tau_i}\Sigma = \Delta H_i - \int_{\tau_{i-1}}^{\tau_{i+1}} \nabla_u H(\psi, x, u) \nabla_{\tau_i} u \, dt$$

and just $\Delta H_i$ for all other nodes. As we now have interior arcs of two types, the rules of generation of 4.2 are appropriately completed with rules of selection between such arcs. The penalty coefficient $\rho$ is gradually increased up to $10^7$ in an outer loop of the algorithm described in section 3.3.

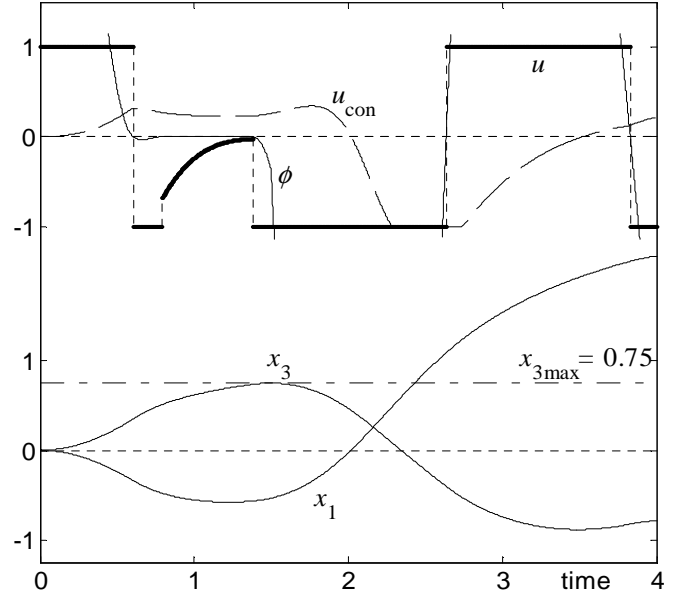

Fig. 3.

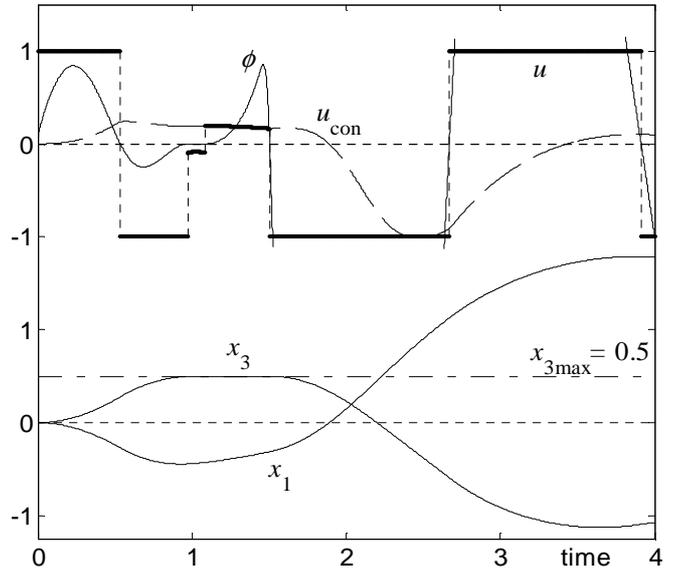

Fig. 4.

The computational experiment is performed for three values of the state bound $x_{3\max}$: 0.75, 0.5 and 0.25. Also in this example the MSE method exhibits rapid convergence. For the sake of brevity we omit the history of structural evolution and present only the final results, see Figs. 3, 4 and 5. For $x_{3\max} = 0.75$ the trajectory of $x_3$ exhibits a touch point, and the optimal control has a considerable singularity and no constrained arc. For a stricter bound $x_{3\max} = 0.5$, the touch point develops into a constrained arc, neighboring with a short singular arc. For $x_{3\max} = 0.25$ there is no singular arc at all, but the constrained arc becomes longer. Every figure shows also the plot of the scaled switching function, here defined as $\phi = -c\psi_2 + \psi_4$. The presented solutions apparently satisfy the optimality conditions of the Maximum Principle.

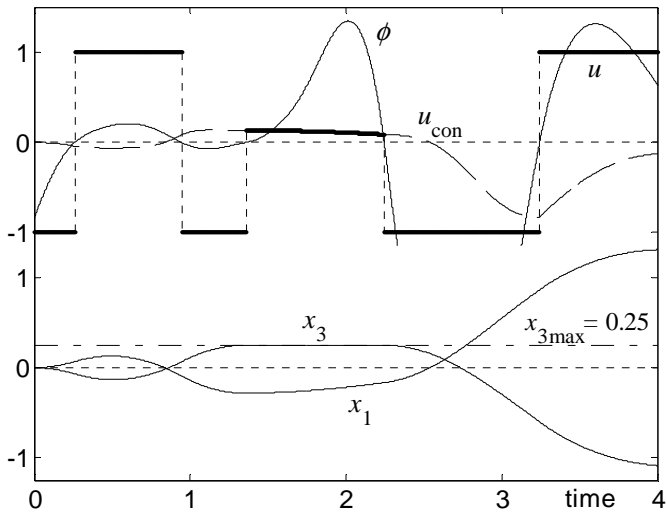

Fig. 5.

## 5 Conclusions

Two features of MSE, i.e., the continuous reconstruction of decision space and increase of its dimension during optimization are also encountered in the well established numerical methods of optimal control. The latter is used in direct algorithms in the form of mesh refinement while the former is a foundation of the homotopy approach. The novelty of MSE lies in the systematic way the available information on expected structure of optimal control is used, and in the concept of control preservation by structural changes. Thanks to these, monotone descent of the performance index is guaranteed and the number of decision variables is economically determined, thus high rates of convergence can be achieved together with relatively accurate representation of optimal control. All intermediate results are control admissible solutions of the control problem (without discontinuous approximations of state trajectories which appear in the homotopy methods with multiple shooting). Besides, among other direct methods, MSE offers the rare possibility of straightforward verification of the maximum principle optimality conditions in case of successful termination.

The computational experience confirms that the rate of convergence of MSE is better than in comparable direct methods, and similar to the rate of convergence of indirect algorithms such as multiple shooting or indirect collocation. On the other hand, the area of convergence of MSE is like in other direct algorithms, and much larger than in indirect methods. The MSE method can be easily extended to problems with free horizon, including time-optimal. It is also possible to treat state constraints in a direct way (without penalty functions), similarly as control constraints and singularities in the presented examples.

## References

[1] J. T. Betts: *Survey of numerical methods for trajectory optimization*. J. Guidance, Control and Dynamics, **21**, 193-207, (1998).

[2] A. E. Bryson, Jr.: *Dynamic Optimization*. Addison – Wesley – Longman, Menlo Park, (1999).

[3] R. Bulirsch, F. Montrone, H. J. Pesch: *Abort landing in the presence of a windshear as a minimax optimal control problem, part 2: multiple shooting and homotopy*. J. Optimization Theory and Applications, Vol. 70, 223- 254, (1991).

[4] J. Kierzenka, L.F. Shampine: *A BVP solver based on residual control and the Matlab PSE*. ACM Transactions on Mathematical Software, **27**, 3, 299-316, (2001).

[5] A. Korytowski, M. Szymkat, A. Turnau: *Time-optimal control of a pendulum on a cart*. CCATIE, Kraków 1998, 177-270 (in Polish).

[6] M. Pauluk, A. Korytowski, A. Turnau, M. Szymkat: *Time optimal control of 3D crane*. Proc. 7th IEEE MMAR 2001, Międzyzdroje, Poland, 28-31 August 2001, 927-932.

[7] I. Queinnec, B. Dahhou: *Optimization and control of a fedbatch fermentation process*. Optimal Control Applications and Methods, Vol. 15, 175-191 (1994).

[8] H. Shen, P. Tsiotras: *Time-optimal control of axi-symmetric rigid spacecraft with two controls*. J. Guidance, Control and Dynamics, **22**, 682-694, (1999).

[9] H.R. Sirisena: *A gradient method for computing optimal bang-bang control*. Int. J. Control, **19**, 257-264, (1974).

[10] O. von Stryk: *User's guide for DIRCOL - a direct collocation method for the numerical solution of optimal control problems*. Ver. 2.1, Technical University of Munich, (1999).

[11] M. Szymkat, A. Korytowski, A. Turnau: *Computation of time optimal controls by gradient matching*. Proc. 1999 IEEE CACSD, Kohala Coast, Hawai'i, August 22-27, 1999, 363-368.

[12] M. Szymkat, A. Korytowski, A. Turnau: *Variable control parameterization for time-optimal problems*. Proc. 8th IFAC CACSD 2000, Salford, U.K., 11-13 September 2000, T4A.

[13] M. Szymkat, A. Korytowski, A. Turnau: *Variable parameterization method for optimal control using second order approximations*. SIAM OP02, Toronto, Canada, 20-22 May 2002, CP 34, 68.

[14] M. Szymkat, A. Korytowski, A. Turnau: *Extended variable parameterization method for optimal control*, Proc. IEEE CCA/CCASD 2002, Glasgow, Scotland, 18-20 September 2002.

[15] A. Turnau, A. Korytowski, M. Szymkat: *Time optimal control for the pendulum-cart system in real-time*. Proc. 1999 IEEE CCA, Kohala Coast, Hawai'i, August 22-27, 1999, 1249-1254.

[16] J. Wen, A.A. Desrochers: *An algorithm for obtaining bang-bang control laws*. J. Dynamic Systems, Measurement and Control, **109**, 171-175, (1987).

[17] http://aq.ia.agh.edu.pl/vpmoc